\documentclass[11pt]{article}
\title{ Further Developments in  Finite Fibonomial Calculus}
\author{Ewa Krot \\
\\Institute of Computer Science, Bia{\l}ystok University\\
PL-15-887 Bia{\l}ystok, ul.Sosnowa 64, POLAND\\
e-mail: ewakrot@wp.pl}

\usepackage{amsmath,amsthm}

\chardef\bslash=`\\ 
\newtheorem{defn}{Definition}[section]
\newtheorem{obs}{Observation}[section]

\newtheorem{rem}{Remark}[section]
\newtheorem{thm}{Theorem}[section]

\begin{document}
\maketitle
\begin{abstract}
Primary definitions, notation and general observations of Finite
Fibonomial Operator Calculus (FFOC) are presented.
A.K.Kwa\'sniewski's combinatorial interpretation of Fibonomial
coefficients by the use of Fibonacci cobweb poset is given. Some
elements of Incidence Algebra of Fibonacci cobweb poset ${\bf P}$
are defined.
\end{abstract}
\section{Introduction}
Finite Fibonomial Operator Calculus (FFOC) is a special case of
Extented Finite Operator Calculus (FOC). The "Calculus of
Sequences" was started by Ward in 1936 \cite{1} and considered by
many authors after him. The very foundations of  FOC were given by
O.V.Viskov, G. Markowsky and finally by the author of \cite{8,2,3}
. The idea of FFOC due to A.K.Kwa\'sniewski \cite{2} (see Example
2.1) is now continued  \cite{4} by the present author.
\section{Primary Definitions, notation and general observations of FFOC}
FFOC ows on the famous Fibonacci sequence
$F=\left\{F_{n}\right\}_{n \geq 1}$ introduced by Leonardo
Fibonacci [Pisano] in 1202 ("Liber Abaci") defined by the
reccurence relation:
\begin{equation}
 \Bigg\{\begin{array}{c}F_{0}=0,\;F_{1}=F_{2}=1\\
F_{n+2}=F_{n+1}+F_{n} \;\;\;\;n \geq 1\end{array}
\end{equation}
or equivalently by the formula (Binet form):
\begin{equation}
F_{n}=\frac{1}{\sqrt{5}}\Bigg(\frac{1+\sqrt{5}}{2}\Bigg)^{n}-\frac{1}{\sqrt{5}}\Bigg(\frac{1-\sqrt{5}}{2}\Bigg)^{n},\;\;\;n
\geq 0.
\end{equation}
Let us define for above sequence  what
 follows:
\renewcommand{\labelenumi}{(\arabic{enumi})}
\begin{enumerate}
\item $F$-factorial:
$$n_{F}!\equiv F_{n}!\equiv F_{n}F_{n-1}...F_{2}F_{1},\;\;\;F_{0}!=1.$$
\item
$$n^{\underline{k}}_{F}\equiv n_{F}\,(n-1)_{F}\,(n-2)_{F}\ldots (n-k+1)_{F}\equiv F_{n}F_{n-1}F_{n-2}\ldots F_{n-k+1}$$
 \item and finally $F$-binomial (Fibonomial ) coefficients:
 \begin{equation}
\binom{n}{k}_{F}=\frac{F_{n}!}{F_{k}!\,F_{n-k}!}=\frac{n^{\underline{k}}_{F}}{k_{F}!},\;\;\;\;\;\;\;\binom{n}{0}_{F}=1.
\end{equation}
\end{enumerate}
In the following we shall consider   the algebra {\bf P} of
polynomials over the field {\bf K} of characteristic zero.

 The
main object of FFOC is the $F$-derivative defined below.
\begin{defn}
 The $F$-derivative is the linear operator   $\partial_{F}:{\bf P}\rightarrow{\bf P}$   such that
   $\partial_{F}x^{n}=F_{n}x^{n-1}$  for $n \geq 0$.
\end{defn}
\begin{defn}
The $F$-translation operator is the linear operator \\
$E^{y}(\partial_{F}): {\bf P}\rightarrow{\bf P}$ of the form:
$$E^{y}(\partial_{F})=\exp_{F}\{ y \partial_{F} \}=
\sum_{k \geq 0} \frac{y^{k} \partial_{F}^{k}}{F_{k}!},\;\;\;\;\;\;
y\in{\bf K}$$ and
$$\forall_{p\in{\bf P}}\;\;\;\;p(x+_{F}y)=E^{y}(\partial_{F})p(x)\;\;\;\;
x,y\in{\bf K}.$$
\end{defn}
 \begin{defn}
A linear operator $T:{\bf P}\rightarrow{\bf P}$ is said to be
 $\partial_{F}$-shift invariant iff
$$\forall_{y \in {\bf K}}\;\;\;\; \left[T,E^{y}(\partial_{F})\right]=
TE^{y}(\partial_{F})-E^{y}(\partial_{F})T=0$$ We shall denote by
$\Sigma_{F}$ the algebra of $F$-linear $\partial_{F}$-shift
invariant operators.
\end{defn}

 \begin{defn}
Let $Q(\partial_{F})$ be a formal series in powers of
$\partial_{F}$ and \\$Q(\partial_{F}):{\bf P}\rightarrow{\bf P}$.
$Q(\partial_{F})$ is said to be $\partial_{F}$-delta operator iff
\renewcommand{\labelenumi}{(\alph{enumi})}
\begin{enumerate}
\item $Q(\partial_{F}) \in \Sigma_{F}$ \item $Q(\partial_{F})(x) =
const \neq 0$.
 \end{enumerate}
\end{defn}
One can show that every $\partial_{F}$-delta operator reduces
degree of any polynomial by one. For every $\partial_{F}$-delta
operator $Q(\partial_{F})$ there exists the uniquely determined
$\partial_{F}$-basic polynomial sequence $\{q_{n}(x)\}_{n\geq 0}$,
($deg\;q_{n}(x)=n$),satysfying conditions:
\renewcommand{\labelenumi}{(\arabic{enumi})}
\begin{enumerate}
\item $q_{0}(x)=1;$ \item $q_{n}(0)=0,\; n\geq 1;$ \item
$Q(\partial_{F})q_{n}(x)=F_{n}q_{n-1}(x),\;\;n \geq 0.$
\end{enumerate}

\begin{thm} {\em (First Expansion Theorem)}\\
Let $T \in \Sigma_{F}$ and let $Q(\partial_{F})$ be a
$\partial_{F}$-delta operator with $\partial_{F}$-basic polynomial
sequence $\{q_{n}\}_{n \geq 0}$. Then
$$T = \sum _{n \geq 0} \frac{a_{n}}{F_{n}!}Q(\partial_{F})^{n};
\quad a_{n} = [Tq_{k}(x)]_{x=0}.$$
\end{thm}

\begin{thm}{\em (Isomorphism Theorem)}\\
Let $\Phi_{F}={\bf K}_{F}[[t]]$ be the algebra of formal {\em
exp}$_{F}$ series in $t \in {\bf K}$ ,i.e.:
$$f_{F}(t)\in \Phi_{F}\;\;\;\;  iff\;\;\;\; f_{F}(t)= \sum_{k \geq 0}
\frac{a_{k}t^{k}}{F_{k}!}\;\;\; for\;\;\; a_{k}\in {\bf K},$$ and
let the $Q(\partial_{F})$ be a $\partial_{F}$-delta operator. Then
$\Sigma_{F} \approx \Phi_{F}$.
 The isomorphism \\ $\phi : \Phi_{F} \rightarrow \Sigma_{F}$ is given by the
 natural correspondence:
 $$f_{F}  \left( {t} \right) = \sum_{k \geq 0}
\frac{a_{k}t^{k}}{F_{k}!}\; \buildrel {into} \over \longrightarrow
\; T_{\partial_{F}} = \sum_{k \geq 0}
\frac{a_{k}}{F_{k}!}Q(\partial_{F})^{k}. $$
\end{thm}
\begin{rem} {\em
In the algebra $\Phi_{F}$ the product is given by the fibonomial
convolution, i.e.:
$$
\left( {\;\sum_{k \geq 0} {\frac{{a_{k}} }{{F_{k}  !}}} x^{k}\;}
\right) \left( {\;\sum_{k \geq 0} {\frac{{b_{k}} }{{F_{k}  !}}}
x^{k}\;} \right)= \left( {\;\sum_{k \geq 0} {\frac{{c_{k}}
}{{F_{k}  !}}} x^{k}\;} \right)$$ where
$$c_{k} = \sum_{l \geq 0} \binom{k}{l}_{F}  a_{l} b_{k -l}.$$ }
\end{rem}
\begin{defn}
A polynomial sequence $\{s_{n}\}_{n\geq 0}$ is called the sequence
of Sheffer $F$-polynomials of the $\partial_{F}$-delta operator
$Q(\partial_{F})$ iff

\renewcommand{\labelenumi}{\em (\arabic{enumi})}
\begin{enumerate}
\item $s_0(x)=const\neq 0$; \item
$Q(\partial_{F})s_{n}(x)=F_{n}s_{n-1}(x);\; n\geq 0.$
\end{enumerate}
\end{defn}
It can be shown that $\{s_{n}\}_{n\geq 0}$ is the sequence of
Sheffer $F$-polynomials of $\partial_{F}$-delta operator
$Q(\partial_{F})$ with $\partial_{F}$-basic polynomial sequence
$\{q_{n}\}_{n  \geq 0}$ iff there exists an invertible $S \in
\Sigma _{F}$ such that $s_{n}(x)=S^{-1}q_{n}(x)$ for  $n\geq 0$.
We shall refer to a given labeled by $\partial_{F}$-shift
invariant invertible operator $S$   Sheffer $F$-polynomial
sequence $\{s_{n}\} _{n\geq 0}$ as the sequence of Sheffer
$F$-polynomials of the $\partial_{F}$-delta operator
$Q(\partial_{F})$ relative to $S$.
\section{Some examples of $F$-polynomials}
The Sheffer F-polynomials form a wide class of new polynomial
sequences. Some examples of them are given below:
\renewcommand{\labelenumi}{\em (\arabic{enumi})}
\begin{enumerate}
\item $\partial_{F}$-basic polynomials
 of the operator $\Delta_{F}=E^{1}(\partial_{F})-I$:
$$q_{n}(x)=
\frac{F_{n}!}{n!}\hat{x}_{F}(E^{-1}(\partial_{F})\hat{x})^{n-1}[1]$$
 where $\hat{1}_{F}^{-1}\hat{x}_{F}=\hat{x}$ for $(\hat{x}p)(x)=xp(x),\;\;
 p\in {\bf P}$.
 \item $\partial_{F}$-basic polynomials of the operator
 $\nabla_{F}=I-E^{-1}(\partial_{F})$:
 $$q_{n}(x)=\frac{F_{n}!}{n!}\hat{x}_{F}(E(\partial_{F})\hat{x})^{n-1}[1].$$
\item $\partial_{F}$-basic polynomials of $F$-Abel operator
$A(\partial_{F})=\partial_{F}E^{a}(\partial_{F})= \sum\limits_{k
\geq 0} \frac{a^{k}}{F_{k}!}\partial_{F}^{k+1}$: $$
A^{(a)}_{n,F}(x)=\frac{F_{n}}{n}\sum_{k=0}^{n-1}
  \binom{n-1}{k}_{F}(-an)^{k}\frac{n-k}{F_{n-k}}x^{n-k}.$$
  Here are the first examples of Abel $F$-polynomials of order $a$:\\
\textrm{}\\
$A^{(a)}_{0,F}(x)=1\\ \textrm{}\\
A^{(a)}_{1,F}(x)=x\\ \textrm{}\\
A^{(a)}_{2,F}(x)=x^{2}+ax\\ \textrm{}\\
A^{(a)}_{3,F}(x)=x^{3}-4ax^{2}+6a^{2}x\\ \textrm{}\\
A^{(a)}_{4,F}(x)=x^{4}-9ax^{3}+48a^{2}x^{2}+48a^{3}x\\ \textrm{}\\
A^{(a)}_{5,F}(x)=x^{5}-20ax^{4}+225a^{2}x^{3}-750a^{3}x^{2}+625a^{4}x\\
\textrm{}\\
A^{(a)}_{6,F}(x)=x^{6}-40ax^{5}+960a^{2}x^{4}-6480a^{3}x^{3}+17280a^{4}x^{2}-
\\
-10368a^{5}x\\ \textrm{}\\
A^{(a)}_{7,F}(x)=x^{7}-78ax^{6}+3640a^{2}x^{5}-50960a^{3}x^{4}+267540a^{4}
x^{3}-499408a^{5}x^{2}+\\+218491a^{6}x\\ \textrm{}\\
A^{(a)}_{8,F}(x)=x^{8}-147ax^{7}+131004a^{2}x^{6}-349440a^{3}x^{5}
+3727360a^{4}x^{4}-\\-13418496a^{5}x^{3}+17891328a^{6}x^{2}-5505024a^{7}x.
\\ \textrm{}\\$
 \item  $\partial_{F}$-basic polynomials of
   $F$-Laguerre operator $L(\partial_{F})=\frac{\partial_{F}}
{\partial_{F}-I}=\sum\limits_{k \geq 0}\partial_{F}^{k+1}$:
   $$   L_{n,F}(x)=\frac{F_{n}}{n}\sum_{k=1}^{n}(-1)^{k}
  \frac{F_{n}!}{F_{k}!}
  \binom{n-1}{k-1}_{F}\frac{k}{F_{k}}x^{k}.$$
  Here are the few first of them:
\\ \textrm{}\\
$L_{0,F}(x)=1\\ \textrm{ } \\
L_{1,F}(x)=-x\\ \textrm{ } \\
L_{2,F}(x)=x^{2}-\frac{1}{2}x\\ \textrm{ } \\
L_{3,F}(x)=-x^{3}+\frac{8}{3}x^{2}-\frac{4}{3}x\\ \textrm{ } \\
L_{4,F}(x)=x^{4}-\frac{27}{4}x^{3}+18x^{2}-\frac{9}{2}x\\  \textrm{ } \\
L_{5,F}(x)=-x^{5}+20x^{4}-135x^{3}+180x^{2}-30x\\  \textrm{ } \\
L_{6,F}(x)=x^{6}-\frac{160}{3}x^{5}+\frac{3200}{3}x^{4}-3600x^{3}+3200x^{2}-
320x\\  \textrm{ } \\
L_{7,F}(x)=-x^{7}+\frac{13182}{7}x^{6}-\frac{54080}{7}x^{5}+
\frac{540800}{7}x^{4}-\frac{1216800}{7}x^{3}+\frac{648960}{7}x^{2}-\\
 \textrm{}\\\;\;\;\;\;\;-\frac{40560}{7}x\\ \textrm{}\\
L_{8,F}(x)=x^{8}-\frac{3087}{8}x^{7}+\frac{223587}{4}x^{6}-1498580x^{5}+
9937200x^{4}-13415220x^{3}+\\ \;\;\;\;+4471740x^{2}-171990x.$
  \item Hermite $F$-polynomials are Sheffer $F$-polynomials of the\\
$\partial_{F}$ -delta operator $\partial_{F}$ relative to
invertible $S \in \Sigma_{F}$ of the form \\ $S=\exp_{F} \{
\frac{a \partial_{F}^{2}}{2}\}$:
$$H_{n,F}(x)=S^{-1}x^{n}=\sum\limits_{k \geq 0}\frac{(-a)^{k}}{2^{k}F_{k}!}
n^{\underline{2k}}_{F}x^{n-2k}.$$
\item  Let $S=(1-\partial_{F})^{-\alpha- 1}$. The Sheffer $F$-polynomials of\\
 $\partial_{F}$-delta operator $L(\partial_{F})=\frac{\partial_{F}}
 {\partial_{F}-1}$  relative to $S$ are Laguerre $F$-polynomials of order
 $\alpha$ are given by formula:
 $$L^{(\alpha)}_{n,F}(x)=\sum\limits_{k \geq 0}\frac{F_{n}!}{F_{k}!}\binom
{\alpha+n}{n-k}_{F}(-x)^{k}$$ for $\alpha \neq -1$. Here are the
few first of Laguerre $F$-polynomials of order $\alpha=1$:
\\ \textrm{}\\
$L^{(1)}_{0,F}(x)=1\\ \textrm{ } \\
L^{(1)}_{1,F}(x)=-x+1\\ \textrm{ } \\
L^{(1)}_{2,F}(x)=x^{2}-2x+2\\ \textrm{ } \\
L^{(1)}_{3,F}(x)=-2x^{3}+6x^{2}-12x+6\\ \textrm{ } \\
L^{(1)}_{4,F}(x)=6x^{4}-30x^{3}+90x^{2}-90x+30\\  \textrm{ } \\
L^{(1)}_{5,F}(x)=-30x^{5}+240x^{4}-1200x^{3}+1800x^{2}-1200x+240\\
 \textrm{ } \\
L^{(1)}_{6,F}(x)=240x^{6}-3120x^{5}+24960x^{4}-62400x^{3}+62400x^{2}-
24960x+3120.$
 \item
Bernoullie's $F$-polynomials of order 1 are Sheffer
$F$-polynomials of
\\$\partial_{F}$ -delta operator $\partial_{F}$ related to
invertible $S=\left(\frac
{\exp_{F}\{\partial_{F}\}-I}{\partial_{F}}\right)^{-1}$. They are
given by formula:
 $$ B_{n,F}(x)= \sum_{k \geq
 0}\frac{1}{F_{k+1}}\binom{n}{k}_{F}x^{n-k}.$$
   Here we give few first of them:\\
\textrm{}\\
$B_{0,F}(x)=1\\ \textrm{ } \\
B_{1,F}(x)=x+1\\ \textrm{ } \\
B_{2,F}(x)=x^{2}+x+\frac{1}{2}\\ \textrm{ } \\
B_{3,F}(x)=x^{3}+2x^{2}+x+\frac{1}{3}\\ \textrm{ } \\
B_{4,F}(x)=x^{4}+3x^{3}+3x^{2}+x+\frac{1}{5}\\  \textrm{ } \\
B_{5,F}(x)=x^{5}+5x^{4}+\frac{15}{2}x^{3}+5x^{2}+x+\frac{1}{8}\\ \textrm{ } \\
B_{6,F}(x)=x^{6}+8x^{5}+20x^{4}+20x^{3}+8x^{2}+x+\frac{1}{13}\\  \textrm{ } \\
B_{7,F}(x)=x^{7}+13x^{6}+52x^{5}+\frac{260}{3}x^{4}+52x^{3}+13x^{2}+x+
\frac{1}{21}\\ \textrm{}\\
B_{8,F}(x)=x^{8}+21x^{7}+\frac{273}{2}x^{6}+364x^{5}+364x^{4}+\frac{273}{2}
x^{3}+21x^{2}+x+\frac{1}{36}\\ \textrm{}\\
B_{9,F}(x)=x^{9}+34x^{8}+357x^{7}+1547x^{6}+\frac{12376}{5}x^{5}+1547x^{4}+
357x^{3}+\\ \textrm{}\\+34x^{2}+x+\frac{1}{55}.$
\end{enumerate}
\begin{rem}{\em
It is known (from the general case, see \cite{8,7,4}) that every
sequence of Sheffer $F$-polynomials is orthogonal in respect to
special inner product associated with this sequence. The problem
of classical orthogonality of these polynomials is under
investigation.}
\end{rem}

\newpage

 \section{Combinatorial Content of FFOC}
 Let us quote some standard interpretations:
 \begin{itemize}
 \item binomial coefficients $\binom{n}{k}$ denote the number
 of $k$-element subsets of  set with $n$ elements;
 \item Stirling numbers of the first kind $\Big[ \begin{array}{c}n\\k
 \end{array}\Big]$denote the number of permutations of $n$
 elements containing exactly $k$ cycles;
 \item Stirling numbers of the second kind
 $\Big\{\begin{array}{c}n\\k\end{array}\Big\}$ denote the number
 of partitions of $n$ elements into $k$ blocks;
 \item q-Gaussian coefficients $\binom{n}{k}_{q}$ denote the
 number of $k$ dimensional subspaces in $n$-th dimensional space
 over Galois field $GF(q)$.
 \end{itemize}
The analogues of above numbers are  fibonomial coefficients
defined in Section 2:
\begin{equation}
\binom{n}{k}_{F}=\frac{F_{n}!}{F_{k}!\,F_{n-k}!}=
\frac{n^{\underline{k}}_{F}}{k_{F}!},\;\;\;\;\;\;\;\binom{n}{0}_{F}=1.
\end{equation}
One can show (using some properties of Fibonacci numbers) that
they are of integer type. Our question is: {\bf Does the
combinatorial (classic like)interpretation of them exist?} It was
unknown until now. The combinatorial interpretation of fibonomial
coefficients has been given by A.K.Kwa\'sniewski in \cite{5,6} by
the use of special partially ordered set (Fibonacci cobweb poset).
His idea was to start with a famous rabbits grown Fibonacci tree.
Let us define a partially ordered set ${\bf P}$ via its Hasse
diagram. It looks like rabbits tree , i.e. it consists of levels
labeled by Fibonacci numbers (the $n$-th level consist of $F_{n}$
elements). Every element of $n$-th level ($n \geq 1, n \in {\bf
N}$) is in partial order relation with every element of the
$(n+1)$-th level but it's not with any element from the level in
which he lies ($n$-th level) except from himself. We will call
this poset ${\bf P}$ the Fibonacci cobweb poset.

We can also define it using the characteristic function $\zeta$ of
the partial order relation $\leq$ in ${\bf P}$. It is defined for
any poset as follows:
$$ \zeta (x,y)=\Big\{\begin{array}{l}1\;\;for\;\;\; x \leq y\\0\;\;otherwise\end{array}$$
The $\zeta$ function of ${\bf P}$ may be expressed by the formula:
$$\zeta(x,y)=\left\{\begin{array}{l}
0\;\;\;if \;x>y\\
0\;\;\;(\sum_{i=1}^{k-1}F_{i})+1 \leq x,y \leq
\sum_{i=1}^{k}F_{k};\;\;x \neq y;\;k\geq 2\\
1\;\;\;x=y\\
1\;\;\;otherwise\end{array}\right.$$ In the above the condition
$(\sum_{i=1}^{k-1}F_{i})+1 \leq x,y \leq \sum_{i=1}^{k}F_{k};\;\;x
\neq y;\;k\geq 2$ means that $x,y$ are the elements of the same
($k$-th) level. We can use the following property of Fibonacci
numbers:
\begin{equation}\label{r5}
\sum_{i=1}^{k}F_{k}=F_{k+2}-1;\;\;\; k\geq 1
\end{equation}
 to get better looking
formula for $\zeta$ function:
\begin{equation}\label{r6}
 \zeta(x,y)=\left\{\begin{array}{l}
0\;\;\;if\; x>y\\0\;\;\;F_{k+1} \leq x,y \leq
F_{k+2}-1;\;\;x \neq y;\;k\geq 3\;\;\;\;\;\;\;\;\;\;\;\;\quad\\
1\;\;\;x=y\\
1\;\;\;otherwise\end{array}\right.
\end{equation}

The incidence $\zeta$ function infinite matrix representing
uniquely cobweb poset ${\bf P}$ was given at first by
A.K.Kwa\'sniewski in \cite{5}. It is of the form
$$\left[ \begin{array}{cccccccccccccccccccccc}
1&1&1&1&1&1&1&1&1&1&1&1&1&1&1&1&1&1&1&1&1&\ldots \\
&1&1&1&1&1&1&1&1&1&1&1&1&1&1&1&1&1&1&1&1&\ldots \\
&&1&0&1&1&1&1&1&1&1&1&1&1&1&1&1&1&1&1&1&\ldots \\
&&&1&1&1&1&1&1&1&1&1&1&1&1&1&1&1&1&1&1&\ldots \\
&&&&1&0&0&1&1&1&1&1&1&1&1&1&1&1&1&1&1&\ldots \\
&&&&&1&0&1&1&1&1&1&1&1&1&1&1&1&1&1&1&\ldots \\
&&&&& &1&1&1&1&1&1&1&1&1&1&1&1&1&1&1&\ldots \\
&&&&&&&1&0&0&0&0&1&1&1&1&1&1&1&1&1&\ldots \\
&&&&&&&&1&0&0&0&1&1&1&1&1&1&1&1&1&\ldots \\
&&&&&&&&&1&0&0&1&1&1&1&1&1&1&1&1&\ldots \\
&&z&e&r&o&s &&&&1&0&1&1&1&1&1&1&1&1&1&\ldots \\
&&&&&&&&&&&1&1&1&1&1&1&1&1&1&1&\ldots \\
&&&&&&&&&&&&1&0&0&0&0&0&0&0&1&\ldots \\
&&&&&&&&&&&&&1&0&0&0&0&0&0&1&\ldots \\
&&&&&&&&&&&&&&1&0&0&0&0&0&1&\ldots \\
&&&&&&&&&&&&&&&1&0&0&0&0&1&\ldots \\
&&&&&&&&&&&&&&&&1&0&0&0&1&\ldots \\
&&&&&&&&&&&&&&&&&1&0&0&1&\ldots \\
.&.&.&.&.&.&.&.&.&.&.&.&.&.&.&.&.&.&.&.&.&\ldots
\end{array}\right]$$

From the definition of ${\bf P}$ follows:
\begin{obs}{\em \cite{5,6}}\\
The number of maximal chains starting from the root (the first
level), to reach $n$-th level (labeled by $F_{n}$) is equal to
$n_{F}!$.
\end{obs}
\begin{obs}{\em \cite{5,6}}\\
The number of maximal chains starting from the any fixed point
from the level labeled by $F_{k}$ to reach any point at the $n$-th
level is equal to\\
$n^{\underline{n-k}}_{F}=n^{\underline{m}}_{F},\;\;(n=k+m)$.
\end{obs}
Let us consider the following question, ( quotation from \cite{5,6}):\\
 {\em "Let us
denote by $P_{m}$ a subposet of $P$ consisting of points up to
$m$-th level points: $\bigcup_{s=1}^{n}\Phi_{s}$; $\Phi_{s}$ is
the set of elements of the $s$-th level. Consider now the
following behaviour of a "sub-cob" moving from any given point of
the $F_{k}$ level of the poset up. It behaves as it has been born
right there and can reach at first $F_{2}$ points up, then $F_{3}$
points up , $F_{4}$ and so on - thus climbing up to to the level
$F_{k+m}=F_{n}$ of the poset ${\bf P}$. It can see- its Great
Ancestor at the root $F_{1}$-th level and potentially follow one
of its own accessible finite subposet $P_{m}$. One of many
$P_{m}$'s rooted at the $k$-th level may be found. How many?"}
\begin{obs}{\em \cite{5,6}}\\
Let $n=k+m$. The number of subposets $P_{m}$ rooted at any fixed
point at the level labeled by $F_{k}$ and ending at the $n$-th
level (labeled by $F_{n}$) is equal to:
$$ \binom{n}{m}_{F}=\binom{n}{m}_{f}=\frac{n^{\underline{k}}_{F}}{k_{F}!}.$$
\end{obs}
\section{The Incidence Algebra ${\bf I(P)}$}
One can define the incidence algebra of ${\bf P}$ (locally finite
partially ordered set) as follows (see \cite{9,10}):
$$ {\bf I(P)}=\{f:{\bf P}\times {\bf P}\longrightarrow {\bf R};\;\;\;\; f(x,y)=0\;\;\;unless\;\;\; x\leq y\}.$$
The sum of two such functions $f$ and $g$ and multiplication by
scalar are defined as usual. The product $H=f\ast g$ is defined as
follows:
$$ h(x,y)=(f\ast g)(x,y)=\sum_{z\in {\bf P}:\;x\leq z\leq y} f(x,z)\ast g(z,y).$$
It is immediately verified that this is an associative algebra
over the real field ( associative ring).

The incidence algebra has an identity element $\delta (x,y)$, the
Kronecker delta. The zeta function of ${\bf P}$ is an element of
${\bf I(P)}$. It was expressed by $\delta$ in \cite{6} from where
we learn that:
$$\zeta =\zeta_{1}-\zeta_{0}$$
where for $x,y \in {\bf N}$:
\begin{equation}
\zeta_{1}(x,y)=\sum_{k=0}^{\infty}\delta(x+k,y)
\end{equation}
\begin{equation}
\zeta_{0}(x,y)=\sum_{k \geq 0}\sum_{s \geq 0}\delta
(x,F_{s+1}+k)\sum_{r=1}^{F_{s}-k-1}\delta (k+F_{s+1}+r,y).
\end{equation}

The knowledge of $\zeta$ enables us to construct other typical
elements of incidence algebra perfectly suitable for calculating
number of chains, of maximal chains etc. in finite sub-posets of
${\bf P}$. The one of them is M\"{o}bius function indispensable in
numerous inversion type formulas of countless applications. It is
known that the $\zeta$ function of a locally finite partially
ordered set is invertible in incidence algebra and its inversion
is so called the M\"{o}bius function $\mu$ i.e.:
$$\zeta \ast \mu=\mu \ast \zeta=\delta.$$The M\"{o}bius function
$\mu$ of Fibonacci cobweb poset ${\bf P}$ is presented here for
the first time and it's due to the present author. It can be
recovered using the following reccurence formula (see \cite{9}):
\begin{equation}
 \left\{\begin{array}{l}\mu(x,x)=1\;\;\;\;for\;\; all \;\;x\in{\bf
 P}\quad \quad \quad \\\ \\\begin{displaystyle}
\mu (x,y)=-\sum_{x\leq z<y}\mu
(x,z)\end{displaystyle}\end{array}\right.
\end{equation}
Then we get:
\begin{equation}
\mu(x,y)=\left\{\begin{array}{l} 0\;\;\;\;\;\;x>y\\
1\;\;\;\;\;\;x=y\\
0\;\;\;\;\;\;\;F_{k+1}\leq x,y\leq F_{k+2}-1;\;x\neq y;\;k\geq 3\\
\\-1\;\;\;\;\;\;F_{k+1}\leq x\leq F_{k+2}-1<F_{k+2}\leq y\leq
F_{k+3}-1\\ \\
-\prod_{l=k+1}^{n-1}(1-F_{l})\quad \quad \quad F_{k+1}\leq x\leq F_{k+2}-1,\\
\quad \quad \quad \quad \quad \quad \quad \quad \quad \quad
F_{n+1}\leq y\leq F_{n+2}-1;\;n>k+1,\end{array}\right.
\end{equation}
where:
\begin{itemize}
\item the condition $F_{k+1}\leq x,y\leq F_{k+2}-1;\;x\neq
y;\;k\geq 3$ means that $x,y$ are different elements of $k$-th
level; \item the condition $F_{k+1}\leq x\leq
F_{k+2}-1<F_{k+2}\leq y\leq F_{k+3}-1$ means that $x$ is an
element of $k$-th level and $y$ is an element of $(k+1)$-th level;
\item the condition $F_{k+1}\leq x\leq F_{k+2}-1,\;F_{n+1}\leq
y\leq F_{n+2}-1;\;n>k+1$ means that $x$ is an element of $k$-th
level and $y$ is an element of $n$-th level.
\end{itemize}

So the $\mu$ function matrix has the following structure:
$$\left[ \begin{array}{cccccccccccccccccc}
1&-1&0&0&0&0&0&0&0&0&0&0&0&0&0&0&0&\ldots \\
&1&-1&-1&1&1&1&-2&-2&-2&-2&-2&8&8&8&8&8&\ldots \\
&&1&0&-1&-1&-1&2&2&2&2&2&-8&-8&-8&-8&-8&\ldots \\
&&&1&-1&-1&-1&2&2&2&2&2&-8&-8&-8&-8&-8&\ldots \\
&&&&1&0&0&-1&-1&-1&-1&-1&4&4&4&4&4&\ldots \\
&&&&&1&0&-1&-1&-1&-1&-1&4&4&4&4&4&\ldots \\
&&&&&&1&-1&-1&-1&-1&-1&4&4&4&4&4&\ldots \\
&&&&&&&1&0&0&0&0&-1&-1&-1&-1&-1&\ldots \\
&&&&&&&&1&0&0&0&-1&-1&-1&-1&-1&\ldots \\
&&&&&&&&&1&0&0&-1&-1&-1&-1&-1&\ldots \\
&&z&e&r&o&s &&&&1&0&-1&-1&-1&-1&-1&\ldots \\
&&&&&&&&&&&1&-1&-1&-1&-1&-1&\ldots \\
&&&&&&&&&&&&1&0&0&0&0&\ldots \\
&&&&&&&&&&&&&1&0&0&0&\ldots \\
&&&&&&&&&&&&&&1&0&0&\ldots \\
&&&&&&&&&&&&&&&1&0&\ldots \\
&&&&&&&&&&&&&&&&1&\ldots \\
&&&&&&&&&&&&&&&&&\ldots \\
.&.&.&.&.&.&.&.&.&.&.&.&.&.&.&.&.&\ldots

\end{array}\right]$$

{\bf {Acknowledgements}}

I would like to thank to Prof. A.K.Kwa\'sniewski for his remarks
and
 guideness.

 AMS Classification numbers: 11C08, 11B37, 47B47

\end{document}